\documentclass[twoside,11pt]{article}

\usepackage{blindtext}
\usepackage[preprint]{jmlr2e}

\usepackage{amssymb}
\usepackage{amsmath}

\usepackage{natbib}
\usepackage{url} 
\usepackage{xcolor}
\usepackage{listings}
\usepackage{amsmath,amsfonts,amssymb}
\usepackage{multirow}
\usepackage{makecell}
\usepackage{units}

\newcommand{\bb}{{\boldsymbol{b}}}
\newcommand{\bx}{{\boldsymbol{x}}}
\allowdisplaybreaks

\lstdefinestyle{mystyle}{
    commentstyle=\color{vscodegreen},
    keywordstyle=\color{blue},
    stringstyle=\color{purple},
    basicstyle=\ttfamily\footnotesize,
    breakatwhitespace=false,         
    breaklines=true,                 
    captionpos=b,                    
    keepspaces=true,                 
    numbersep=5pt,                  
    showspaces=false,                
    showstringspaces=false,
    showtabs=false,                  
    tabsize=2
}

\usepackage{hyperref}

\usepackage{lastpage}

\ShortHeadings{}{}
\firstpageno{1}

\begin{document}

\title{Improving  hp-Variational Physics-Informed Neural Networks for Steady-State Convection-Dominated Problems}

\author{\name Thivin Anandh $^{1,\dagger}$ \email thivinanandh@iisc.ac.in \\
       \name Divij Ghose $^{1,\dagger}$ \email divijghose@iisc.ac.in \\
        \name Himanshu Jain $^{1}$\email himanshuj1@iisc.ac.in \\
        \name Pratham Sunkad $^{2}$\email prathamsunkad@smail.iitm.ac.in \\
        \name Sashikumaar Ganesan$^{1}$ \email sashi@iisc.ac.in \\
        \name Volker John $^{3,4}$ \email john@wias-berlin.de \\
       \addr $^1$ Department of Computational and Data Sciences\\
       Indian Institute of Science, Bangalore\\
       Karnataka, India \\
      \addr $^2$ Department of Mechanical Engineering\\
       Indian Institute of Technology, Madras\\
       Tamil Nadu, India \\
      \addr $^3$ Numerical Mathematics and Scientific Computing\\
       Weierstrass Institute for Applied Analysis and Stochastics\\
       Berlin, Germany \\
       \addr $^4$ Department of Mathematics and Computer Science\\
       Freie University\\
       Berlin, Germany \\
       $^{\dagger}$ denotes equal contribution}

\editor{My editor}

\maketitle

\begin{abstract}
This paper proposes and studies two extensions of applying hp-variational physics-informed 
neural networks, more precisely the FastVPINNs framework, to convection-dominated convection-diffusion-reaction problems. 
First, a term in the spirit of a SUPG stabilization is included in the loss functional and a
network architecture is proposed that predicts spatially varying stabilization parameters.
Having observed that the selection of the indicator function in hard-constrained Dirichlet 
boundary conditions has a big impact on the accuracy of the computed solutions, the second 
novelty is the proposal of a network architecture that learns good parameters for a
class of indicator functions. Numerical studies show that both proposals 
lead to noticeably more accurate results than approaches that can be found in the literature.

\end{abstract}


\section{Introduction}
Convection-diffusion-reaction (CDR) problems are fundamental models for simulating transport events. 
CDR problems capture the interaction between a fluid's bulk motion, or convection, its progressive spreading of properties by random 
molecular motion, or diffusion, and the impact from other quantities in coupled systems, which might be modelled as 
the reaction term. They constitute a framework for modelling the transport of variables like temperature or concentration.

Let $\Omega \subset \mathbb{R}^2$ be a bounded domain with polygonal Lipschitz-continuous boundary $\partial\Omega$. The Lebesgue and Sobolev spaces on this domain are denoted by $L^p(\Omega)$ and $W^{k,p}(\Omega)$, respectively, where, $1\leq p \leq \infty$, $k \geq 0$. The Hilbert space, equivalent to $W^{k,2}(\Omega)$, is denoted by $H^k(\Omega)$. 
Then, a linear CDR boundary value problem, already in nondimensionalized form,  is given by
\begin{align}
    \begin{split}
        -\varepsilon\Delta u(\bx) + \bb(\bx)\cdot\nabla u(\bx) + c(\bx)u(\bx) &= f(\bx), \quad \text{in} \  \Omega, \\
        u(\bx) &=  g(\bx), \quad\text{on} \  \partial\Omega.
    \end{split}
    \label{eq:cdr_eqn}
\end{align}   
Here, $\bx =(x,y) \in \overline\Omega$, $u(\bx)$ is the unknown scalar solution and $f(\bx) \in L^2(\Omega)$ is a known source function. In addition, $\varepsilon \in \mathbb{R}^+$ is the diffusion 
coefficient, $\bb \in (W^{1, \infty}(\Omega))^2$ is the convection field, and $c \in L^{\infty}(\Omega)$ is the reaction field. The Dirichlet boundary 
condition is prescribed by  $g(\bx) \in H^{1/2}(\partial\Omega)$. 

Of particular interest in many applications is the situation that the convective term dominates the equation. In particular, it is often 
larger than the diffusive term by several orders of magnitude, i.e., the CDR problem is said to be convection-dominated when
$\varepsilon \ll L\|\bb\|_{L^\infty(\Omega)}$, where $L$ is the characteristic length of the problem. With the Péclet number, $\mathrm{Pe}$, defined as 
\begin{equation*}
    \mathrm{Pe} = \frac{L\|\bb\|_{L^\infty(\Omega)}}{\varepsilon},
\end{equation*}
the convection-dominated regime is often characterized by $\mathrm{Pe} \gg 1$, in many applications $\mathrm{Pe} \gtrsim 10^6$.
In the convection-dominated situation, which is also sometimes called singularly 
perturbed, typical solutions of CDR problems possess so-called layers, which are thin regions where the solution has very 
steep gradients, e.g., see \cite{roos2008robust}. It is known from asymptotic analysis that the layer width is $\mathcal O(\varepsilon)$
or $\mathcal O(\sqrt{\varepsilon})$, depending on the type of the layer.

Classical discretization techniques for CDR problems, based on finite element, finite volume, or finite difference methods, 
face in the convection-dominated regime the difficulty that the layer width is (much) smaller than the affordable mesh 
width. Consequently, the layers, which are the most important features of the solution, cannot be resolved. The situation 
that important features of a solution cannot be resolved is typical for multiscale problems. For convection-dominated 
CDR problems, the layers are subgrid scales. It is well known that one has to utilize so-called stabilized discretizations 
to cope with this difficulty. In the framework of finite element methods, one of the earliest but still one of the most 
popular proposals is the Streamline-Upwind Petrov--Galerkin (SUPG) method from \cite{HB79,BH82}. The numerical analysis of
this linear discretization is well understood, see \cite{roos2008robust}. 
However, it requires the choice of a function of stabilization parameters, which is often chosen to be a piecewise constant 
function on the given triangulation. For these parameters, usually proposals based on a one-dimensional problem with constant 
coefficients are applied, e.g., see \cite{john2007spurious}. However, numerical solutions with such parameters often exhibit
spurious oscillations in a vicinity of layers. It is shown in \cite{john2011posteriori,JKW23} that (in some sense) optimized 
stabilization parameters might reduce the size of  spurious oscillations considerably. A more recent approach for determining 
appropriate stabilization parameters is presented in \cite{yadav2023artificial}, which utilizes artificial neural networks. 
Both techniques, the optimization and the use of neural networks, are inherently nonlinear. The recent review \cite{BJK24}
studied methods that lead to numerical solutions which satisfy discrete maximum principles, i.e., which provide numerical 
solutions without unphysical values. It is emphasized in this paper that to achieve this property in combination with high accuracy, then one has 
to apply a nonlinear method, because different approaches have to be applied for the subgrid and for the large scales. 
As conclusion, convection-dominated CDR problems are a challenging class of problems from the numerical point of view
and the accurate and physically consistent numerical solution requires to use nonlinear methods. 

In recent times, deep learning based methods are increasingly being used in applied mathematics problems, 
such as solving boundary value problems with partial differential equations (PDEs). Collectively termed scientific machine learning (SciML), 
\cite{cuomo2022scientific, osti_1478744, psaros2023uncertainty}, such methods act in conjunction with or entirely 
replace classical approaches. 
The field of SciML has witnessed a surge in development, leading to a growing number of accessible SciML libraries, e.g., \cite{lu2021deepxde, modulus}.
The application of such methods has grown exponentially since the introduction of physics-informed neural networks (PINNs) in \cite{712178, raissi2019physics}. 
In addition to the typical data-driven loss functional used in neural networks, a PINN incorporates an additional loss term to minimize the residual of 
the underlying PDE and to enforce the physics-based constraints. The abilities to obtain gradients conveniently 
with automatic differentiation  and to train the same neural network for both forward and inverse modelling are some characteristics 
that make PINNs an interesting alternative of classical methods, \cite{abueidda2021meshless, lu2021physics}.

Variational PINNs~(VPINNs) are an extension of PINNs that use a weak form of the PDE in the loss functional, see 
\cite{kharazmi2019variational, khodayi2020varnet}. VPINNs use, as PINNs, a neural network for approximating the solution. 
The test functions belong to a polynomial function space, which is analogous to the Petrov--Galerkin framework 
in finite element methods. Usually, the domain $\Omega$ is decomposed into subdomains (elements) and the support of each test 
function is just one subdomain, defining a so-called hp-VPINN. 
Similar to finite element methods, the accuracy of hp-VPINNs has been shown, for certain problems, 
to increase by further domain decomposition (h-refinement) or by increasing the order of the polynomial space (p-refinement), see \cite{kharazmi2021hp}.

However, the standard use of hp-VPINNs shows some bottlenecks and limitations. A bottleneck is the high computational expense of the training, e.g., 
as reported in \cite{FHJ24}, which increases in direct proportion to the number of elements. Furthermore, many current hp-VPINNs implementations face 
difficulties with complex geometries,  which include skewed quadrilateral elements, which are prevalent in real-world scenarios. 
To overcome these issues, FastVPINNs, as introduced in~\cite{anandh2024fastvpinns, anandh2024fastvpinns-joss,ghose2024fastvpinns}, employs a tensor-based method for loss functional computations, 
coupled with the implementation on GPUs, 
leading to a remarkable reduction in training time of order $\mathcal O(100)$ compared to conventional hp-VPINNs. Additionally, FastVPINNs use bilinear transformations,  
allowing to simulate problems on complex geometries with skewed quadrilateral elements.

Despite all the positive aspects of PINNs and hp-VPINNs, it has been observed in \cite{krishnapriyanCharacterizingPossibleFailure2021} that 
these approaches face difficulties when applied to singularly perturbed problems, like convection-dominated CDR problems. Many contributions that study PINNs or (hp-)VPINNs 
for CD(R) problems consider either the (much more) simpler one-dimensional case or only a mildly convection-dominated regime, 
see the introduction of \cite{FHJ24} for a survey of corresponding works. The choice of examples for the current paper is guided by 
\cite{FHJ24,FJZ23,MJZ24}, where strongly convection-dominated problems are studied. The emphasis of these papers was on using 
loss functionals that are alternative to the standard residual loss, the selection of collocation points in PINNs, 
and on approaches that preserve given bounds of the solution. For two 
examples studied in \cite{FHJ24,FJZ23} the obtained results were characterized as nonsatisfactory. 

The first contribution of this paper consists in proposing alternative loss functionals, compared with the standard 
residual loss, for convection-dominated problems. A term in the spirit of the SUPG stabilization is introduced, 
which contains a stabilization parameter. In addition, the effect of including a regularization with respect to
the weights of the neural network is studied. A first main novelty is the proposal of a FastVPINNs architecture 
that learns spatially varying stabilization parameters. It is shown that this approach leads to noticeable 
improvements of the accuracy of the numerical solutions. The second main contribution is on the realization of
hard-constrained Dirichlet boundary conditions, more precisely, on the definition of the indicator function 
that separates the extension of the boundary conditions and the function that is learnt by the network. 
Also for this task, a FastVPINNs  architecture is proposed that learns the parameters for a class of indicator 
functions. The adaptively chosen parameters lead to a considerable increase of the accuracy 
compared with manually chosen parameters. 

The paper is organized as follows. Section~\ref{sec:setup} describes the general setup of the FastVPINNs framework. 
The loss functionals are introduced in Section~\ref{sec:loss}. Numerical studies with constant stabilization 
parameters and with manually chosen indicator functions are presented in Section~\ref{sec:num_const_param}.
The two main novelties, the FastVPINNs  architectures for spatially varying stabilization parameters 
and for the adaptive choice of the parameters in the indicator function, are introduced and studied in 
Section~\ref{sec:new_exp}. Finally, Section~\ref{sec:conclusion} summarizes the main results of this paper.

\section{Setup of the FastVPINNs Framework}\label{sec:setup}

This section starts by introducing the variational form of the CDR problem \eqref{eq:cdr_eqn}.
Then, the fast hp-VPINNs are described. 

Let
\[
 V:=\left\{v \in H^1(\Omega)\ : \ v = 0 \ \text{on} \ \partial \Omega \right\},
\]
then the variational form is obtained in the usual way by multiplying equation~\eqref{eq:cdr_eqn} with 
an arbitrary test function $v\in V$, integrating over $\Omega$, and then utilizing integration by parts on the second order 
derivative term, e.g.,  see \cite{ganesan2017finite}. Given an extension $u_{\mathrm{ext}} \in H^1(\Omega)$ of the boundary
condition into $\Omega$.
Then, the variational form of \eqref{eq:cdr_eqn} reads as follows:
Find $u \in H^1(\Omega)$ such that $u-u_{\mathrm{ext}} \in V$ and 
\[
a(u,v) = f(v) \quad \forall\ v \in V, 
\]
where the bilinear form $a(\cdot, \cdot)\:\ H^1(\Omega)\times V \rightarrow \mathbb{R}$ is defined as 
\begin{align}
    \begin{split}
    a(u,v) &:= \int_{\Omega} \varepsilon\nabla u \cdot \nabla v \ \mathrm{d}\bx + \int_{\Omega} \bb \cdot \nabla u \ v \ \mathrm{d}\bx + 
    \int_{\Omega} cu \ v \ \mathrm{d}\bx, \\
        f(v) &:= \int_{\Omega} f \ v \ \mathrm{d}\bx.
    \end{split}
    \label{eqn:cdr_weakform}
\end{align}

For defining a hp-VPINN, the domain $\Omega$ is then divided into an array of non-overlapping cells, labeled as $K_k$, 
where $k=1,2,\ldots, N_{\mathrm{elem}}$, such that their union 
$\bigcup_{k=1}^{N_{\mathrm{elem}}} K_k = \overline\Omega$. 
Let $V_h$ be a finite-dimensional subspace of $V$, spanned by $N_{\mathrm{test}}$ basis functions. 
As a result, the discretized form of \eqref{eqn:cdr_weakform} can be written as follows:
Find $u_h$ such that $u_h - u_{\mathrm{ext},h} \in V_h$ and
\begin{equation}
a_h(u_h,v_h) = f_h(v_h) \quad \forall\ v_h \in V_h,  \label{eqn:cdr_disform}
\end{equation}
where $u_{\mathrm{ext},h}$ is an appropriate discrete extension of the Dirichlet boundary condition and 
\begin{align*}
    \begin{split}
    a_h(u_h,v_h) &:= \sum_{k=1}^{N_{\mathrm{elem}}}\int_{K_k} \varepsilon \nabla u_h \cdot \nabla v_h \ \mathrm{d}\bx + \sum_{k=1}^{N_{\mathrm{elem}}}\int_{K_k} \bb \cdot \nabla u_h \ v_h \ \mathrm{d}\bx \\
    &\quad + \sum_{k=1}^{N_{\mathrm{elem}}}\int_{K_k} cu_h \ v_h \ \mathrm{d}\bx, \\
        f_h(v_h) &:= \sum_{k=1}^{N_{\mathrm{elem}}}\int_{K_k} f \ v_h \ \mathrm{d}\bx.
    \end{split}
\end{align*}
These integrals are approximated by employing numerical quadrature: 
\begin{align*}
    \int_{K_k} \varepsilon\nabla u_h \cdot \nabla v_h \; \mathrm{d}\bx &\approx  \sum_{q=1}^{N_{\mathrm{quad}}} w_q \ \varepsilon\nabla u_h(\bx_q) \cdot \nabla v_h(\bx_q)\ ,  \\
    \int_{K_k} \bb \cdot \nabla u_h \ v_h \ \mathrm{d}\bx &\approx  \sum_{q=1}^{N_{\mathrm{quad}}} w_q \ \bb \cdot \nabla u_h(\bx_q) \ v_h(\bx_q)\ , \\
    \int_{K_k} cu_h \ v_h \ \mathrm{d}\bx &\approx \sum_{q=1}^{N_{\mathrm{quad}}} w_q \ cu_h(\bx_q) \ v_h(\bx_q)\ , \\
       \int_{K_k} f\ v_h \ \mathrm{d}\bx & \approx   \sum_{q=1}^{N_{\mathrm{quad}}}w_q \ f(\bx_q)\ v_h(\bx_q)\ .
\end{align*}
Here, $N_{\mathrm{quad}}$ is the number of quadrature points in a cell, $\{\bx_q\}$ is the set of quadrature nodes,  and $\{w_q\}$ is the set of weights (which include the area of $K_k$).

The hp-VPINN framework \citep{kharazmi2021hp} utilizes specific test functions $v_k$, where $k$ ranges from $1$ to $N_{\mathrm{elem}}$, which are localized and defined 
within individual non-overlapping cells 
\begin{equation*}
    v_k= 
    \begin{cases}
      v^p \neq 0, & \text{in $K_k$,} \\
      0, & \text{elsewhere.}
    \end{cases}
\end{equation*}
Here, $v^p$ represents a polynomial function of degree $p$. 
This selection of test and solution spaces resembles a Petrov--Galerkin finite element method.

By utilizing these functions, we define the cell-wise residual of the variational form~\eqref{eqn:cdr_disform} with $u_{\text{NN}}(\bx; \theta_W,\theta_b)$ by
\begin{align}
    \begin{split}
        \mathcal{W}_k(\theta_W,\theta_b) = \sum_{q=1}^{N_{\mathrm{quad}}} &w_q \big[ \varepsilon\nabla u_{\text{NN}}(\bx_{kq}) \cdot \nabla v_k(\bx_{kq}) + \bb \cdot \nabla u_{\text{NN}}(\bx_{kq}) \; v_k(\bx_{kq})  \\
        &+ c\;u_{\text{NN}}(\bx_{kq}) \; v_k(\bx_{kq}) - f(\bx_{kq})\,v_k(\bx_{kq}) \big] \;.
    \end{split}
\label{eqn:cdr_weak_residual}
\end{align}
Now, the variational loss is given by 
\begin{equation} 
    \mathcal{L}_{\mathrm{var}}(\theta_W,\theta_b) = \frac{1}{N_{\mathrm{elem}}}\sum_{k=1}^{N_{\mathrm{elem}}}\left|\mathcal{W}_k(\theta_W,\theta_b)\right|^2.
    \label{eq:var_loss}
\end{equation}

There are two conceptually different ways of incorporating Dirichlet boundary conditions in the framework of hp-VPINNs. 
First, these boundary conditions can be learnt. To this end, a term is included in the loss functional that measures 
the difference between the prescribed boundary values and the values of the neural network solution $u_{\text{NN}}(\bx;\theta_W,\theta_b)$ 
taken at the boundary. In this approach, the Dirichlet boundary conditions are satisfied usually only approximately. 
This is, in our opinion, not appropriate in the case of convection-dominated CDR problems. First, conditions posed on 
the inlet boundary $(\bb\cdot\boldsymbol{n})(\boldsymbol{x}) < 0$, $\boldsymbol{x}\in\partial\Omega$, where $\boldsymbol{n}$ is the outward pointing unit normal vector at 
$\partial\Omega$, are transported in the domain and thus are essential for determining the correct form of the solution inside the domain. And second, 
boundary layers, which are an essential feature of the solution, might not be present in  $u_{\text{NN}}(\bx;\theta_W,\theta_b)$  if the 
boundary condition is only approximately satisfied. In addition, one has to choose a weight for the boundary term that 
relates the importance of this term to other terms in the loss functional. For these reasons, we decided to pursue the 
alternative approach of using a hard-constrained imposition of the boundary conditions in the hp-VPINNs. With this approach, 
the boundary conditions, which are given data of the problem, are satisfied exactly. 

The ansatz function for using hard-constrained boundary conditions, see \cite{lu2021physics}, is 
\begin{equation}
    u^{\text{hard}}_{\text{NN}}(\bx; \theta_W, \theta_b) = j(\bx) + h(\bx)u_{\text{NN}}(\bx; \theta_W, \theta_b),
    \label{eq:hard_ansatz}
\end{equation}
Here, $j(\bx)$ is a continuous function in $\Omega$ that satisfies $j(\bx)=g(\bx)$ for $\bx \in \partial\Omega$ and $h(\bx)$ is an indicator function satisfying 
\begin{equation}\label{eq:indicator_fct}
    h(\bx)= 
    \begin{cases}
      0, & \bx \in \partial\Omega, \\
      > 0, & \text{elsewhere.}
    \end{cases}
\end{equation}
For example, the indicator function used in~\cite{FHJ24} has the form 
\begin{equation}\label{eq:indicator_FHJ24}
    h(\bx) = \left( 1 - \text{e}^{-\kappa x} \right) \left( 1 - \text{e}^{-\kappa y} \right) \left( 1 - \text{e}^{-\kappa (1 - x)} \right)  \left( 1 - \text{e}^{-\kappa (1 - y)} \right),
\end{equation}
where $\kappa$ depends on the  diffusion coefficient $\varepsilon$, but not on the shape of the solution close to the boundary. 
The current paper employs different indicator functions for different examples, 
which are detailed in the respective numerical sections. In addition, we study a method where the parameters for a function of type \eqref{eq:indicator_FHJ24} are obtained from the neural network, in Section~\ref{sec:new_exp}.
Now, the representation \eqref{eq:hard_ansatz} is inserted in \eqref{eqn:cdr_weak_residual} to get the residual $\mathcal{W}^{\text{hard}}_k(\theta_W,\theta_b)$, 
which then can be substituted in \eqref{eq:var_loss} to define the cost functional of the neural network by
\begin{equation}
    \mathcal{L}^{\text{hard}}_{\text{var}}(\theta_W,\theta_b) = \frac{1}{N_{\mathrm{elem}}}\sum_{k=1}^{N_{\mathrm{elem}}}\left|\mathcal{W}^{\text{hard}}_k(\theta_W,\theta_b)\right|^2.
    \label{eq:vpinn_hard_loss}
\end{equation}

\begin{figure}[t!]
    \centering
    \includegraphics[width=0.75\textwidth]{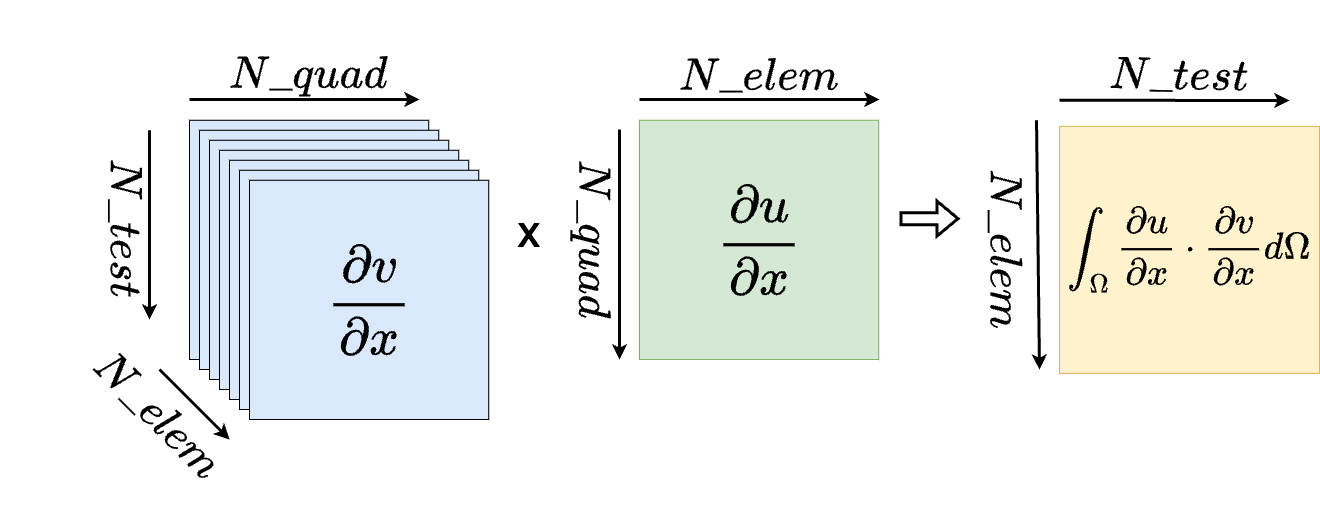}
    \caption{Tensor-based loss computation schematic for FastVPINNs}
    \label{fig:FastVPINNs}
\end{figure}

FastVPINNs computes the variational loss~\eqref{eq:vpinn_hard_loss} by stacking the test functions and their gradients into 
a three-dimensional tensor and reshaping the neural network gradients into a matrix, see Figure~\ref{fig:FastVPINNs}, thus enabling the calculation of the final 
residual vector for all elements in a single operation using BLAS routines available in TensorFlow~\cite{tensorflow2015-whitepaper}. 
This approach has two main advantages: using BLAS operations on the GPU allows for efficient and fast calculations, and stacking 
the gradients and performing tensor-based operations enable the loss for the entire domain to be computed in a single calculation, 
reducing the training time's dependency on the number of elements, which is a major challenge in conventional hp-VPINNs frameworks. 
Additionally, FastVPINNs can handle complex geometries, making this method suitable for a wide range of problems. 
For further details on the implementation of the FastVPINNs framework, see~\cite{anandh2024fastvpinns-joss}.

\section{Loss Functional}\label{sec:loss}
The main topic of \cite{FHJ24} is the investigation of different loss functionals for PINNs and 
hp-VPINNs applied to the numerical solution of convection-dominated CDR problems. The motivation for 
studying different loss functionals stems from the observation that the standard residual is 
not an appropriate choice for performing parameter optimization in stabilized finite element methods, 
see \cite{john2011posteriori}. It was observed in \cite{FHJ24} that it was possible to obtain 
more accurate solutions with loss functionals that are different to the standard residual loss. 

In the current paper, two alternative loss functionals will be studied.

\subsection{SUPG Stabilization Loss}

\begin{figure}[t!]
    \centering
    \includegraphics[width=0.999\textwidth]{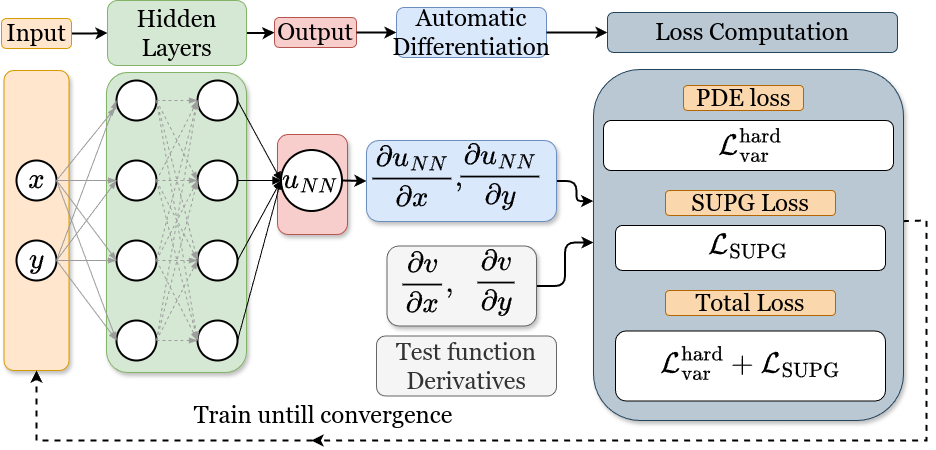}
    \caption{hp-VPINNs architecture for convection-dominated problems with SUPG stabilization, prescribed stabilization parameter, and hard constraints}
    \label{fig:supg_vpinns}
\end{figure}

The SUPG finite element method, which is already mentioned in the introduction, adds an additional 
term to the standard Galerkin finite element discretization. This term essentially introduces 
numerical diffusion in streamline direction. The global SUPG stabilization term has the form 
\begin{equation}
   \mathcal{L_{SUPG}}  = \int_{\Omega} \tau(\bx)\left(\bb\cdot\nabla u(\bx) + cu(\bx)-f(\bx)\right)\left(\bb\cdot\nabla v(\bx)\right)\ \mathrm{d}\bx, 
    \label{eqn:SUPG_Loss}
\end{equation}
where $\tau(\bx)$ is called stabilization parameter and the diffusive term is neglected in the residual (the first 
factor), which is appropriate in the convection-dominated regime. In the framework of finite elements, 
\eqref{eqn:SUPG_Loss} is localized by decomposing the integral in a sum of integrals over the mesh cells 
and then $\tau(\bx)$ is often defined to be a piecewise constant function on the given triangulation. 
In the framework of neural networks, the SUPG stabilization loss functional is given by 
\begin{equation}
    \mathcal{L^{SUPG}_\tau} = \mathcal{L}^{\text{hard}}_{\text{var}} + \mathcal{L_{SUPG}},
    \label{Eq: Supg_loss}
\end{equation}  
where $\mathcal{L}^{\text{hard}}_{\text{var}}$ is the variational loss defined in  \eqref{eq:vpinn_hard_loss}.
Minimizing the loss functional \eqref{Eq: Supg_loss} penalizes large streamline derivatives. 
Hence, the additional SUPG loss can be also interpreted as a physics-informed regularization term, since 
functions with small residuals but large streamline derivatives are excluded from the set of
possible solutions. It is well known that regularization techniques might 
help to prevent overfitting. An architecture of the network for the loss functional including the 
SUPG term with prescribed stabilization parameter is presented in Figure~\ref{fig:supg_vpinns}.

\subsection{$L^2$ Weight Regularization Loss}

The training of a neural network requires the solution of a large-scale non-convex optimization 
problem. For such problems, usually a regularization term is included in the loss functional. 
In addition, this term might counteract overfitting and thus it enhances the generalization capacity of the network.  
A standard approach consists in adding a $L^2$-type regularization, so that the loss functional becomes
\begin{equation}
        \mathcal L_{\lambda}^{\text{reg}} =\mathcal{L}^{\text{hard}}_{\text{var}}  + \frac{\lambda}{N}\sum_jw_j^2
        \label{Eq: Loss_1},
    \end{equation}
where $N$ is the total number of entries in the weight matrices of the network
and $\{w_j\}_{j=1}^N$ is the set of all weights in the network.  
The $L^2$ weight decay regularization parameter $\lambda$ needs to be tuned.

\section{Numerical Studies with Constant Parameters}\label{sec:num_const_param}
In this section, we first validate the FastVPINNs code by solving the Eriksson--Johnson problem, as presented in \cite{sikora2023physics}, 
by comparing accuracy and computing times. 
Subsequently, we address two convection-dominated problems,
assessing the performance of the proposed loss functionals  
\eqref{Eq: Supg_loss} and \eqref{Eq: Loss_1} with constant stabilization and regularization parameter, respectively, within the framework of hp-VPINNs and we compare the results and efficiency 
with those reported 
in the literature.

All simulations were conducted on a system equipped with an AMD Threadripper CPU and an NVIDIA A6000 GPU. For neural network training, we utilized TensorFlow 2.0, specifically employing the \texttt{tf.float64} datatype to ensure high precision. The network architecture incorporates \texttt{tanh} activation functions, and we used the Adam optimizer for the training iterations. We adopted test functions similar to those described in~\cite{kharazmi2021hp}. These functions are defined as $v_k = P_{k+1} - P_{k-1}$, where $P_k$ represents the $k$th order Legendre polynomial. 
For numerical integration, we implemented Gauss--Lobatto--Legendre quadrature routines as used in~\cite{kharazmi2021hp}. 

\subsection{Studied Problems}

\subsubsection*{The Eriksson--Johnson Problem ($\text{P}_{\text{EJ}}$)}\label{sec:appendix_EJ} 

\begin{figure}
    \centering
    \includegraphics[width=\linewidth]{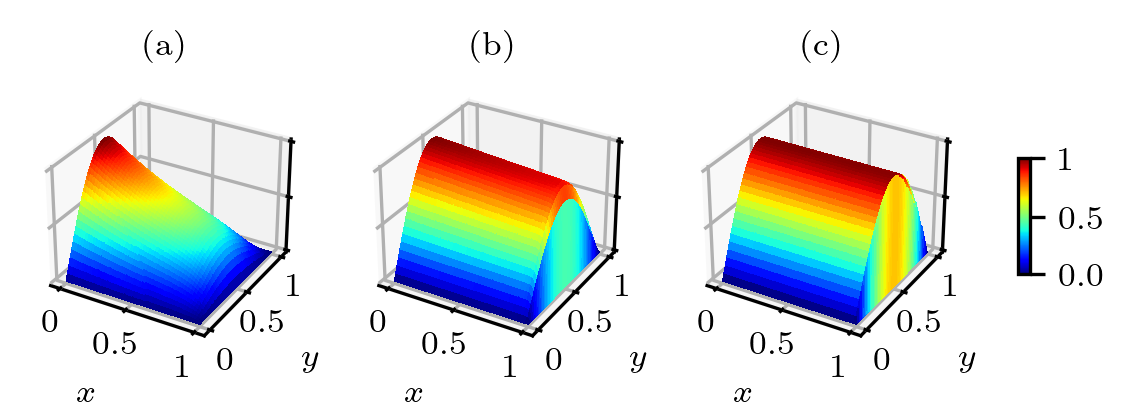}
    \caption{Exact Solution of $\text{P}_{\text{EJ}}$ for (a) $\varepsilon=0.1$ , (b) $\varepsilon=0.01$ and (c) $\varepsilon=0.001$}
    \label{fig:exact_johnson}
\end{figure}

This problem, originally proposed in \cite{EJ91}, was used in \cite{sikora2023physics} for assessing a PINN for 
(mildly) convection-dominated problems. It is defined by $\Omega=(0,1)^2$, $\bb=(1,0)^T$, $c=0$, and the right-hand
side of \eqref{eq:cdr_eqn} is chosen such that the solution has the form 
\begin{align*}
    u(\bx) = \; & \dfrac{\text{e}^{r_1(x-1)}-\text{e}^{r_2(x-1)}}{\text{e}^{-r_1}-\text{e}^{-r_2}}\sin(\pi y) \quad \mbox{with} \\ 
    & r_1 = \dfrac{1+\sqrt{1+4\varepsilon^2\pi^2}}{2\varepsilon},\quad
    r_2 = \dfrac{1-\sqrt{1+4\varepsilon^2\pi^2}}{2\varepsilon}.
\end{align*}
Dirichlet conditions were imposed on $\partial\Omega$. Solutions for various values of the diffusion coefficient 
are displayed in Figure~\ref{fig:exact_johnson}. The smallest value of the diffusion coefficient considered in 
\cite{sikora2023physics}  is $\varepsilon = 0.001$. For small values of the diffusion coefficient, 
the solution exhibits a boundary layer at $x=1$. 

\subsubsection*{Solution with Outflow Boundary Layers ($\text{P}_{\text{out}}$)}\label{sec:outflow}
 
This example, proposed in \cite{JMT97}, features two layers at the outflow boundary, which, according to an asymptotic analysis, have a width of $\mathcal{O}(\varepsilon)$. 
The data for problem \eqref{eq:cdr_eqn} are given by $\Omega = (0,1)^2$, $\varepsilon = 10^{-8}$, $\mathbf{b} = (2,3)^T$, $c = 1$, and the prescribed solution has the form:
\begin{align}
    \begin{split}
        u(x,y) := x y^{2} - y^{2}\exp\left({\frac{2(x-1)}{\varepsilon}}\right) - x\exp\left({\frac{3(y-1)}{\varepsilon}}\right) \\
\qquad\quad + \exp\left({\frac{2(x-1)+3(y-1)}{\varepsilon}}\right).
    \end{split}
    \label{eq:Outflow_layer}
\end{align}
Dirichlet boundary conditions are imposed on $\partial\Omega$. By construction, the outflow boundary layers are situated at $x=1$ and at $y=1$, and there is even a so-called corner singularity at the right upper corner of the domain. Figure~\ref{fig:exact_P}(a) depicts the prescribed solution \eqref{eq:Outflow_layer}.



\subsubsection*{Parabolic layer problem ($\text{P}_{\text{para}}$)}\label{sec:parabolic} 
\begin{figure}
    \centering
    \includegraphics[width=0.99\textwidth]{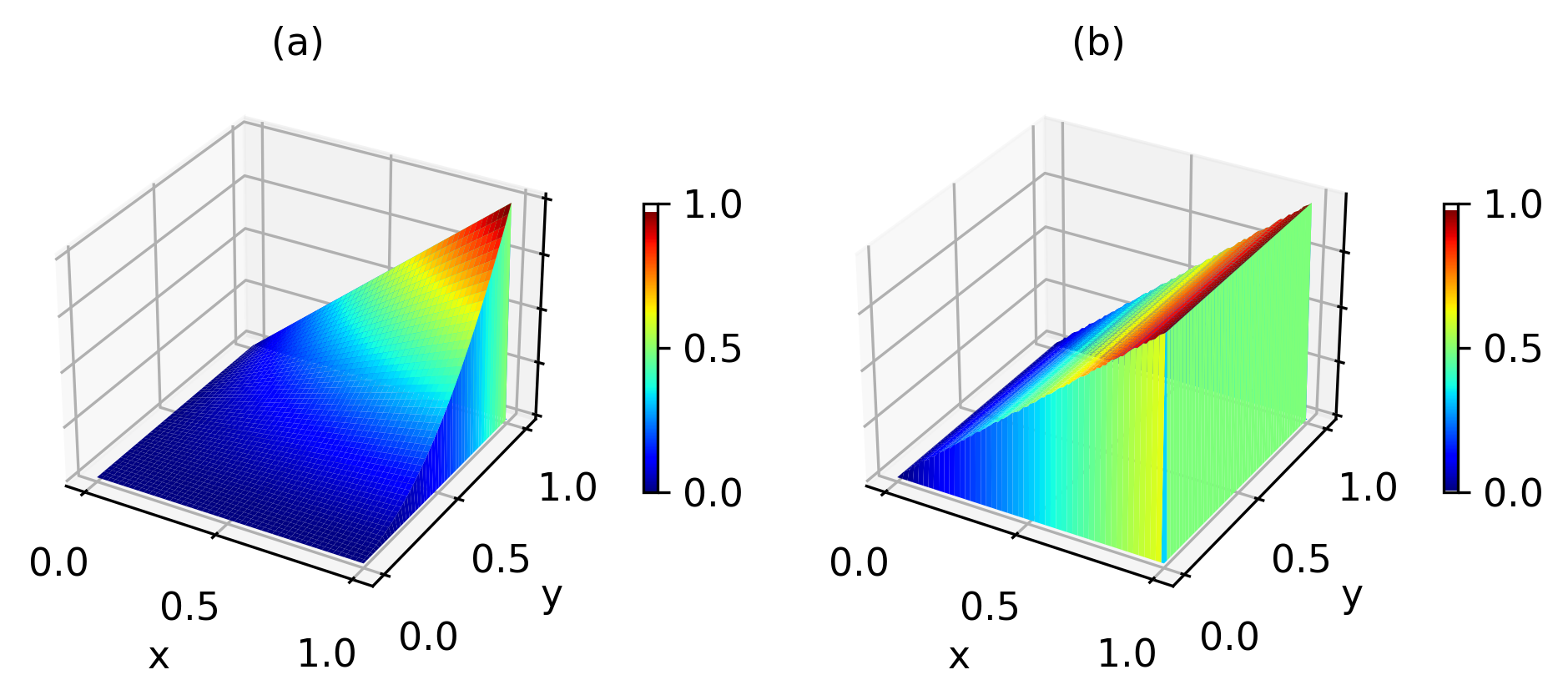}
    \caption{Exact solution for (a) $\text{P}_{\text{out}}$ and (b)  $\text{P}_{\text{para}}$}
    \label{fig:exact_P}
\end{figure}

The data for problem \eqref{eq:cdr_eqn} are given by 
$\Omega = (0,1)^2$, $\varepsilon = 10^{-8}$, $\bb = (1,0)^T$, $c = 0$, and $f=1$. Homogeneous Dirichlet conditions are 
prescribed at $\partial\Omega$. 
The solution exhibits an exponential layer 
at the outflow boundary $x=1$ and two parabolic boundary layers at the characteristic boundaries $y=0$ and $y=1$, as shown in Figure~\ref{fig:exact_P} (b).

To obtain a high accuracy solution for this problem to compare with, we employed the Monotone Upwind-type Algebraically Stabilized (MUAS) 
method. The MUAS method, formulated in~\cite{DBLP:journals/corr/abs-2111-08697}, is an algebraically stabilized scheme. 
It has been shown in the survey \cite{BJK24} that algebraically stabilized scheme are currently the most promising 
finite element approaches for computing solutions of convection-dominated convection-diffusion-reaction problems that satisfy 
Discrete Maximum Principles. 

\subsection{Validation of the FastVPINNs Implementation and Comparison with the Literature}\label{sec:EJ}

In this section, we consider the Eriksson--Johnson problem $\text{P}_{\text{EJ}}$ as studied in~\cite{sikora2023physics} and compare the available 
data with our results. The approach in~\cite{sikora2023physics} utilizes a single element with 6400 quadrature points and 
a weak boundary constraint for the Dirichlet boundaries. In contrast, our hp-VPINNs method employs an $8 \times 8$ element domain with 
$9$ test functions and $100$ quadrature points per element, i.e., in total also 6400 quadrature points, and it uses a hard boundary constraint 
for imposing Dirichlet conditions. Despite these differences, both approaches share common features in their neural network architecture 
and training process. Specifically, both use a neural network with $4$ hidden layers and $20$ neurons per hidden layer, employ the 
Adam optimizer for training with a learning rate of $0.00125$, and run for $40,000$ training epochs.
In this study, we have used $\varepsilon \in \{0.1,0.01,0.001\}$, as it was done in \cite{sikora2023physics}.

The indicator functions for the hard constraints \eqref{eq:hard_ansatz} were chosen to be 
\begin{eqnarray*}
        j(\bx) &= &  \sin(\pi y)\cos\left(\frac{\pi}{2} x\right), \\ 
        h(\bx) &=& \left(1 - \text{e}^{-\kappa_1 x}\right)\left(1 - \text{e}^{-\kappa_1 y}\right)\left(1-\text{e}^{-\kappa_2(1 - x)}\right)\left(1-\text{e}^{-\kappa_1(1 - y)}\right)
\end{eqnarray*}
with $\kappa_1 = 30$ and $\kappa_2 = 10/\varepsilon$. That means, for the indicator function $h(\bx)$ we prescribed a sharp boundary layer at the boundary $x=1$ to account for the 
corresponding layer of the solution. This is not necessary at the other parts of the boundary. The standard loss functional $\mathcal{L}^{\text{hard}}_{\text{var}}$ defined in \eqref{eq:vpinn_hard_loss}, 
which takes into account the hard constraint boundary conditions, was used.


\begin{table}[t!]
\centering
\caption{Eriksson--Johnson problem $\text{P}_{\text{EJ}}$, $L^2_{\text{err}}$
and comparison with the literature} 
\label{table:Error_johnson}
\renewcommand{\arraystretch}{1.2} 
\begin{tabular}{|l|c|l|l|}
\hline
quadrature points & $\varepsilon$ & FastVPINNs & \cite{sikora2023physics} \\  
\hline
6400 & 0.1 & $1.400\cdot10^{-3}$ & 0.265\\
6400 & 0.01 & $1.440\cdot10^{-2}$ & 0.346\\ 
6400 & 0.001 & $7.000\cdot10^{-4}$  & 0.332 \\
\hline
\end{tabular}
\end{table}

\begin{table}[t!]
\caption{Eriksson--Johnson problem $\text{P}_{\text{EJ}}$, average training time per epoch in $\unit{s}$ and comparison with the literature}
\label{tab:EJ_times}
\centering
\begin{tabular}{|c|c|c|}
\hline
FastVPINNs (CPU) & FastVPINNs (GPU) & \cite[Fig.~14]{sikora2023physics} \\
\hline
0.011 & 0.0027 & 0.0375 \\
\hline
\end{tabular}
\end{table}

Tables~\ref{table:Error_johnson} and~\ref{tab:EJ_times} present the obtained results. It can be seen that, on the one hand, the application of FastVPINNs gives 
results that are more accurate by several orders of magnitude than the results from~\cite{sikora2023physics}. Here, the approximation $L^2_{\text{err}}$ of the error in $L^2(\Omega)$ is 
defined by 
\begin{equation}\label{eq:MSE}
L^2_{\text{err}} = \left( \frac{1}{N_T} \sum_{t=1}^{N_T} \left(\left(u-u^{\text{hard}}_{\text{NN}}\right)(\bx_t)\right)^2\right)^{1/2},
\end{equation}
where $N_T$ is the number of test points at which the error is calculated. In all of the examples, the test points $\bx_t$ are given by a $100 \times 100$ grid of 
equidistant points in $\overline\Omega$.
And on the other hand, the training times are considerably smaller, by a factor of about 3 on CPUs and by a factor of 14 if GPUs are used. 

The obtained results validate our implementation of FastVPINNs. In our opinion, the improved accuracy compared with the literature can be attributed to two factors.
First, the use of hard constraints for the boundary conditions, which seems to be essential for boundary layer problems, and second the $h$-refinement approach of hp-VPINNs.

\subsection{Studies with the SUPG Stabilization and the Regularization with Constant Parameters}
\label{subsec:Initial Experimentation}

In this section, we analyze the performance of the SUPG-based loss functional and loss functional with regularization
applied to problems $\text{P}_{\text{out}}$ and $\text{P}_{\text{para}}$. 
Concretely, our studies explore the   
SUPG stabilization $\mathcal{L^{SUPG}_\tau}$ defined in \eqref{Eq: Supg_loss} and PDE loss with regularization $\mathcal{L}_{\lambda}^{\text{reg}}$ given in \eqref{Eq: Loss_1}. For the simulations, we employed a neural network architecture 
consisting of seven hidden layers with 30 neurons each. The domain was discretized into 64 cells, with 100 quadrature points and 36 test functions per cell, 
resulting in a total of 6400 quadrature points. To conduct the experiments, learning rates of $0.01\cdot 3^{-2}$ and $0.01\cdot 3^{-3}$ have been used. The $L^2_{\text{err}}$ is calculated as given in \eqref{eq:MSE}.
The network is trained for 100,000 epochs, and the smallest error obtained during this training, i.e., after each epoch and not necessarily after the final epoch, is reported below. 
The used parameters are similar to those used in the studies described in~\cite{FHJ24}. 
For both problems, we applied an indicator function as presented in~\cite{FHJ24}: 
\begin{equation}
\begin{split}
     h(x,y) &= \left( 1 - \text{e}^{-\kappa x} \right) \left( 1 - \text{e}^{-\kappa y} \right) \left( 1 - \text{e}^{-\kappa (1 - x)} \right)  \left( 1 - \text{e}^{-\kappa (1 - y)} \right), \\
     \kappa &= 10/10^{-8}\;=\;10^9.
\end{split}
\label{eq: Ansatz_initial}
\end{equation}

\begin{table}[t!]
\centering
\caption{Problem $\text{P}_{\text{out}}$, results for using constant parameters in the SUPG loss or regularization loss;
best results reported in \cite{FHJ24}:  $3.419\cdot 10^{-2}$ (PINN) and  $3.619 \cdot 10^{-2}$ (hp-VPINN)}
\label{tab:outflow_layer_params}
\renewcommand{\arraystretch}{1.2}
\small
\begin{tabular}{|l|r|r|r|}
\hline
loss & search range & optimal value & \multicolumn{1}{c|}{$L^2_{\text{err}}$}  \\
\hline
$\mathcal{L}_{\lambda}^{\text{reg}}$ & $[10^{-5}, 10^{-3}]$ & $2.4 \cdot 10^{-4}$ $(\lambda)$ & $3.202\cdot 10^{-2}$ \\
\hline
$\mathcal{L^{SUPG}_\tau}$ & $[10^{-5},10^{-1}]$ & $3.5 \cdot 10^{-2}$ $(\tau)$ & \textbf{$\boldsymbol{1.973\cdot 10^{-2}}$} \\
\hline
\end{tabular}
\end{table}

\begin{table}[t!]
\centering
\caption{Problem $\text{P}_{\text{para}}$, results for using constant parameters in the SUPG loss or regularization loss}
\label{tab:circular_interior_layer_params}
\renewcommand{\arraystretch}{1.2}
\small
\begin{tabular}{|l|r|r|r|}
\hline
loss & search range & optimal value & \multicolumn{1}{c|}{$L^2_{\text{err}}$}
\\
\hline
$\mathcal{L}_{\lambda}^{\text{reg}}$ & $[10^{-8}, 10^{-2}]$ & $4 \cdot 10^{-6}$  $(\lambda)$ & $4.888 \cdot 10^{-1}$ \\
\hline
$\mathcal{L^{SUPG}_\tau}$ & $[10^{-5}, 5]$ & $1.2$ $(\tau)$ & {\textbf{$\boldsymbol{5.384 \cdot 10^{-2}}$}} \\
\hline
\end{tabular}
\end{table}

\begin{figure}[t!]
    \centering
    \includegraphics[width=0.99\textwidth]{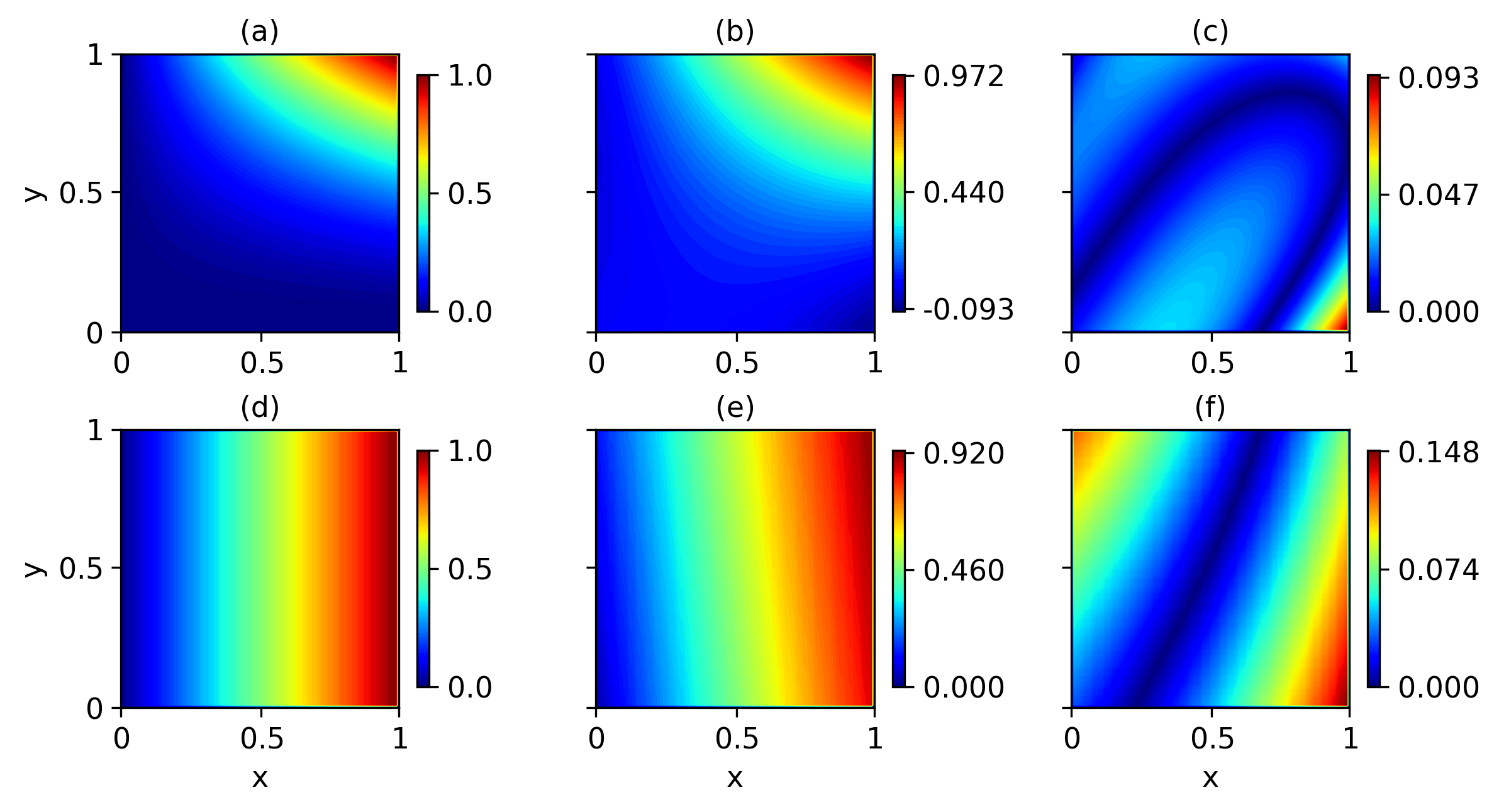}
    \caption{Best results for $\text{P}_{\text{out}}$  (top) and $\text{P}_{\text{para}}$  (bottom). The exact solutions for $\text{P}_{\text{out}}$  and $\text{P}_{\text{para}}$
    are shown in (a) and (d). The predicted solutions are shown in (b) and (e) and the point-wise errors in (c) and (f)}
    \label{fig:Inital_expt}
\end{figure}

Our numerical studies were conducted in two phases. First, we performed an initial broad-range hyperparameter search, computing $L^2_{\text{err}}$ for a wide range of hyperparameters,
as detailed in Tables~\ref{tab:outflow_layer_params} and~\ref{tab:circular_interior_layer_params}. Following this, we conducted a refined hyperparameter 
search based on the best performing values from the initial search. The optimal hyperparameters for each loss functional are presented in aforementioned tables and the 
corresponding solutions and point-wise errors in Figure~\ref{fig:Inital_expt}. 
For both problems, it can be seen that the SUPG stabilization ($\mathcal{L^{SUPG}_\tau}$)
outperformed the approach that uses only regularization ($\mathcal{L}_{\lambda}^{\text{reg}}$). In particular, 
for problem $\text{P}_{\text{out}}$ our methods achieved lower errors compared to the PINNs  and
hp-VPINNs simulations performed in~\cite{FHJ24}. However, on the other hand, both results are not yet completely
satisfactory. There is a comparatively big error in the right lower corner of the domain for the 
numerical solution of $\text{P}_{\text{out}}$  and 
the contour lines of the numerical solution of $\text{P}_{\text{para}}$ are not parallel to the $y$-axis. 
But all in all, we think that the obtained results demonstrate that 
including the SUPG stabilization term is an appropriate way for enhancing the accuracy of variational physics-informed neural networks. 

\begin{figure}
    \centering
    \includegraphics[width=0.8\linewidth]{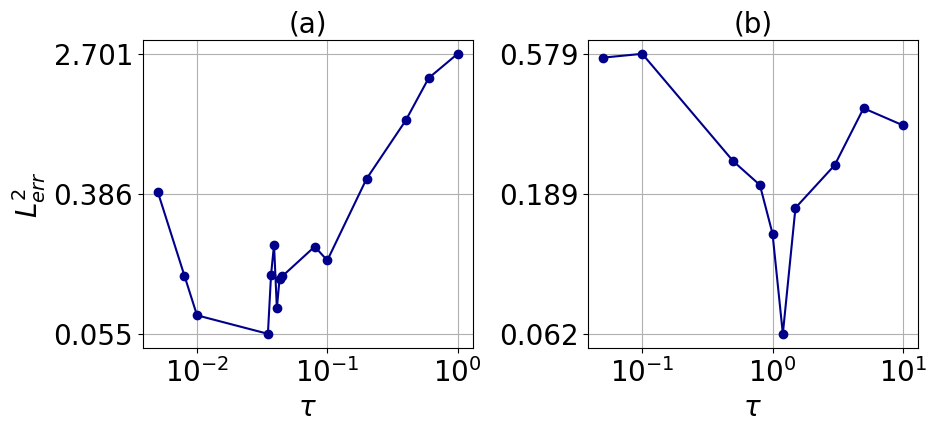}
    \caption{Sensitivity analysis on the $\tau$ parameter for problems $\text{P}_{\text{out}}$ and $\text{P}_{\text{para}}$ using $\mathcal{L^{SUPG}_\tau}$ as loss functional}
    \label{fig:sensitivity}
\end{figure}

Next, results of a sensitivity analysis concerning the accuracy of the numerical solutions for 
variations of the parameter $\tau$ in  $\mathcal{L^{SUPG}_\tau}$ will be presented. 
It is important to note that even when all parameters are kept constant, the minimum error obtained varies noticeably 
across multiple runs of the simulations. To obtain a more representative picture, we repeated the simulations 5 times and used 
the average of the minimum errors to plot the values shown in Figure~\ref{fig:sensitivity}.
It is important to clarify that the errors mentioned in Tables~\ref{tab:outflow_layer_params} and~\ref{tab:circular_interior_layer_params} 
are the result of the best run only.
For $\text{P}_{\text{out}}$, it can be seen that there is a certain, but relatively small, interval for the stabilization parameter
where the obtained errors are close to the best one. In the case of $\text{P}_{\text{para}}$, there is even just 
a single stabilization parameter that gives the by far most accurate result. The best values from Figure~\ref{fig:sensitivity}
correspond with those from Tables~\ref{tab:outflow_layer_params} and~\ref{tab:circular_interior_layer_params}. This study shows 
that a refined search is needed to find the optimal value of the stabilization parameter to obtain accurate results. This search is usually 
not feasible in practice since the solution of the problem is usually not known. 

\section{Learning the Parameters in the SUPG Loss and the Indicator Function}
\label{sec:new_exp}
\begin{figure}[t!]
    \centering
    \includegraphics[width=0.99\textwidth]{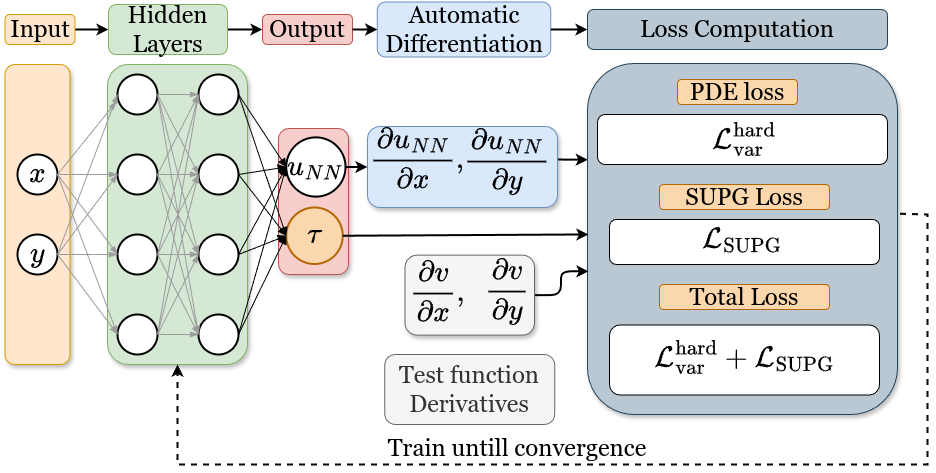}
    \caption{hp-VPINNs architecture for convection-dominated problems with SUPG stabilization that predicts the stabilization parameter $\tau$ and imposes hard constraints}
    \label{fig:supg_tau_vpinns}
\end{figure}

The studies from Section~\ref{sec:num_const_param} and further 
preliminary numerical simulations have identified two main challenges in addressing convection-dominated CDR problems 
with the hp-VPINNs as studied in this paper. The first challenge consists in determining the optimal, or at least a good,
SUPG stabilization parameter $\tau$. The second challenge involves the choice of a suitable indicator function~\eqref{eq:indicator_fct} 
for enforcing Dirichlet boundary conditions, see also \cite{MJZ24} for a discussion of this challenge. 
This second challenge is particularly important because the indicator function 
must not only match the Dirichlet boundary values at the boundary, but its gradient near the boundary should also fit
to the solution's gradient. To address these challenges, we propose two enhancements to hp-VPINNs, and with that in 
particular to FastVPINNs. Firstly, we introduce a modified neural network that simultaneously predicts a spatially varying stabilization 
parameter $\tau(x,y)$ along with the solution as shown in Figure~\ref{fig:supg_tau_vpinns}. Secondly, we also propose an adaptive indicator function whose slope near the boundaries is controlled by a learnable parameter, which is learnt during the training process. 

\subsection{Neural Network Predicting the Stabilization Parameter \texorpdfstring{$\tau(x,y)$}{tau(x,y)}}

The stabilization parameter that controls the artificial streamline diffusion can be expressed as a function 
of the coefficients of the CDR problem and the local mesh width, as it is done in finite element methods. However, the 
difficulty of this approach is that the optimal stabilization parameter is not known for problems defined in 
two- and higher-dimensional domains. In a neural network context, the stabilization parameter $\tau$ can be considered as an additional 
hyperparameter that has to be found by performing hyperparameter tuning, as it was done in Section~\ref{sec:num_const_param}.
However, this approach is time-consuming and it might be infeasible in practice. 

To address this challenge, we propose a neural network that can predict the stabilization parameter $\tau$ along with the solution, 
as shown in Figure~\ref{fig:supg_tau_vpinns}. With this approach, the neural network-predicted $\tau$ can be spatially varying, 
and this variation depends not only on the coefficients of the problem but also on the local shape of the predicted solution. 
This feature can also be observed for optimized SUPG stabilization parameters, e.g., see \cite{john2011posteriori,JKW23}. In this way, 
numerical stabilization is added only where needed.

The network-predicted $\tau$ is constrained using a sigmoid activation function to ensure that it remains positive.
Preliminary numerical studies showed that it is of advantage, for obtaining consistent values of $\tau$ across different 
runs, to set $\tau$ to zero at the boundary of the domain. That means, the quadrature nodes on $\partial\Omega$ do not 
contribute to the SUPG stabilization loss $\mathcal L_{\mathcal{SUPG}}$  defined in \eqref{eqn:SUPG_Loss}. In fact, the
solution on $\partial\Omega$ is known and applying a stabilization there is not mandatory. 
We used the indicator function 
\begin{equation*}
    \tau^{\text{hard}}_{\text{NN}}\left(\bx; \theta_W, \theta_b\right) = w(\bx) \sigma\left(\tau_{\text{NN}}(\bx; \theta_W, \theta_b)\right),
\end{equation*}
where, $\tau_{\text{NN}}(\bx; \theta_W, \theta_b)$ is the stabilization parameter  predicted by the neural network, $\sigma$ is the sigmoid function, 
and $w(\bx)$ is an indicator function defined by 
\begin{equation*}
    w(\bx) = \tanh(50x)  \tanh(50y)   \tanh(50(1-x))  \tanh(50(1-y)).
\end{equation*} 
Then we set 
\begin{equation*}
    \tau(\bx; \theta_W, \theta_b) = \tau_g \tau^{\text{hard}}_{\text{NN}}(\bx; \theta_W, \theta_b),
\end{equation*}
where $\tau_g$, termed as `$\tau$ growth', is controlling the strength of the stabilization factor. The hyperparameter $\tau_g$ 
needs to be fine-tuned. The predicted $\tau$ is then used in the calculation of $\mathcal{L^{SUPG}_\tau}$ given in \eqref{Eq: Supg_loss}. 

We will use problems $\text{P}_{\text{out}}$ and $\text{P}_{\text{para}}$ to demonstrate the effect of applying the learnt  $\tau$ in the SUPG stabilization loss. 
Three approaches will be compared: PDE loss  \eqref{eq:vpinn_hard_loss} only, PDE loss with SUPG stabilization \eqref{Eq: Supg_loss} and constant $\tau$, and PDE  loss with SUPG stabilization \eqref{Eq: Supg_loss} with learnt
$\tau_g$-based SUPG stabilization. 

\begin{figure}[t!]
    \centering
\includegraphics[width=1\textwidth]{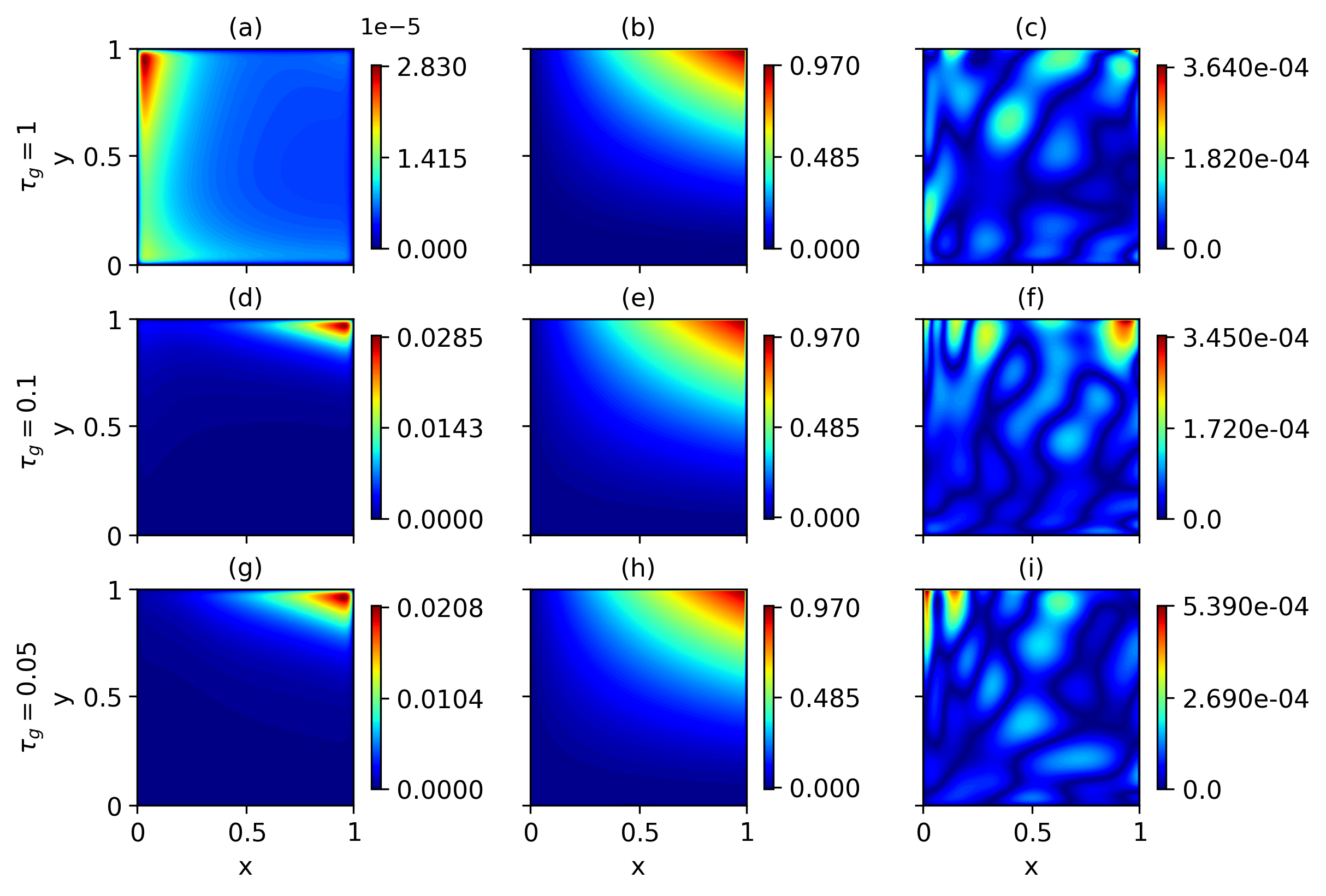}
    \caption{Results for problem $\text{P}_{\text{out}}$  with $\mathcal{L^{SUPG}_\tau}$ for $\tau_g = 1,\; 0.1,\; 0.05$.
   Note that the optimal solution was obtained for $\tau_g = 1$. The stabilization parameter $\tau(x,y)$ is shown in (a), (d) and (g). The predicted solution is depicted in (b), (e) and (h). The point-wise error is presented in (c), (f) and (i)} 
    \label{fig:Modified_ansatz_solution_tau}
\end{figure}

\begin{figure}[h!]
    \centering
\includegraphics[width=1\textwidth]{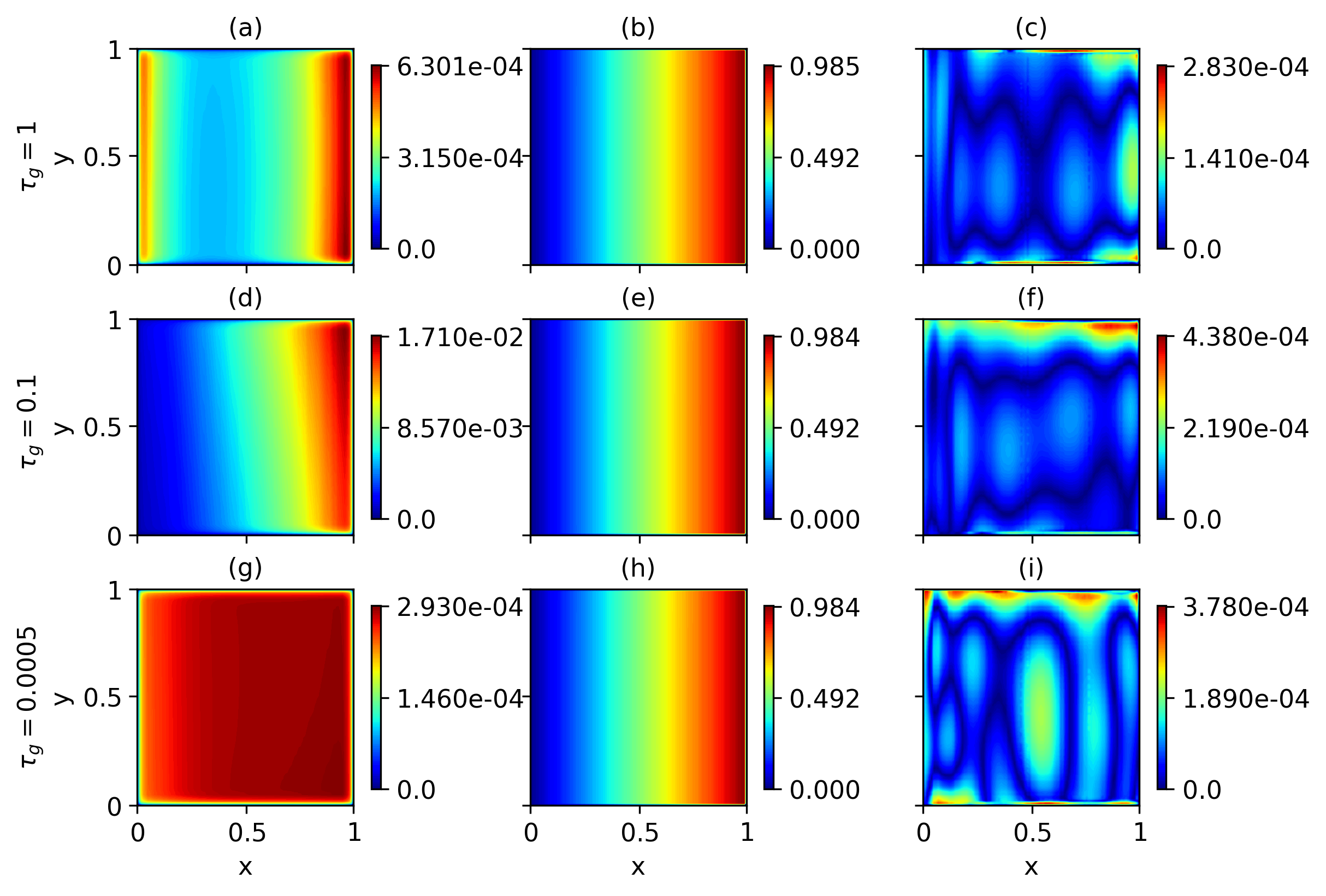}
    \caption{Results for problem $\text{P}_{\text{para}}$  with $\mathcal{L^{SUPG}_\tau}$ for $\tau_g = 1,\; 0.1,\; 0.0005$.
   Note that the optimal solution was obtained for $\tau_g = 1$. The stabilization parameter $\tau(x,y)$ is shown in (a), (d) and (g). The predicted solution is depicted in (b), (e) and (h). The point-wise error is presented in (c), (f) and (i)} 
    \label{fig:sang_Modified_ansatz_solution_tau}
\end{figure}

In our previous studies we observed that using a symmetric indicator function as given in~\eqref{eq: Ansatz_initial} might generate large errors
on the boundary when there is no steep gradient in the solution, compare Figure~\ref{fig:Inital_expt}. Therefore, for this study, we applied a 
modified indicator function that has a steep slope at the outflow boundaries $x=1$ and $y=1$  and a low slope at the other two boundaries:
\begin{equation}
\begin{split}
     h(x,y) &= \left( 1 - \text{e}^{-\kappa_1 x} \right) \left( 1 - \text{e}^{-\kappa_1 y} \right) \left( 1 - \text{e}^{-\kappa_2 (1 - x)} \right)  \left( 1 - \text{e}^{-\kappa_2 (1 - y)} \right), \\
     \kappa_1 &= 30,\;\kappa_2 = 10/10^{-8}\;=\;10^9.
\end{split}
\label{eq: Ansatz_2}
\end{equation}
Similarly for $\text{P}_{\text{para}}$, we applied a modified indicator function, which has a steep slope at boundaries $y=0$, $y=1$ and $x=1$, and a low slope at the boundary $x=0$:
\begin{equation}
\begin{split}
     h(x,y) &= \left( 1 - \text{e}^{-\kappa_1 x} \right) \left( 1 - \text{e}^{-\kappa_2 y} \right) \left( 1 - \text{e}^{-\kappa_2 (1 - x)} \right)  \left( 1 - \text{e}^{-\kappa_2 (1 - y)} \right), \\
     \kappa_1 &= 30,\;\kappa_2 = 10/10^{-8}\;=\;10^9.
\end{split}
\label{eq: Ansatz_3}
\end{equation}
A neural network with seven hidden layers and 30 neurons per layer was used. The network was trained for $50,000$ epochs with a learning rate of $0.01 \cdot 3^{-2}$. 
We used $8 \times 8$ cells with $100$ quadrature points and $36$ test functions per cell. To determine the optimal hyperparameters for the constant $\tau$ and neural network-predicted $\tau$, we initially performed a single run for each hyperparameter within the broad search ranges specified in Tables~\ref{tab:adaptive_tau_accuracy_new_out} and~\ref{tab:adaptive_tau_accuracy_new_para} for the problems $\text{P}_{\text{out}}$ and $\text{P}_{\text{para}}$, respectively.  From these initial runs, we identified the hyperparameters that performed the best for each approach. Subsequently, we ran 10 simulations using these optimal hyperparameters for each study. The average $L^2_{\text{err}}$ values from these 10 runs are presented in Tables~\ref{tab:adaptive_tau_accuracy_new_out} and~\ref{tab:adaptive_tau_accuracy_new_para}.  It can be observed that the model with the learnt $\tau$ (with the hyperparameter $\tau_g$) leads to notably more accurate results
than both the model with constant $\tau$ (as a hyperparameter) and the model without any SUPG loss. Some representative solutions
for the approach with the learnt stabilization parameter are presented in Figures~\ref{fig:Modified_ansatz_solution_tau} and~\ref{fig:sang_Modified_ansatz_solution_tau}.

\begin{table}[t!]
\centering
\caption{$L^2_{\text{err}}$ for problem $\text{P}_{\text{out}}$ using different loss functionals}
\label{tab:adaptive_tau_accuracy_new_out}
\renewcommand{\arraystretch}{1.2}
    \small
    \begin{tabular}{|l|c|c|c|}
    \hline
    loss & $\mathcal{L}^{\text{hard}}_{\text{var}}$ & 
    \begin{tabular}[c]{@{}c@{}}$\mathcal{L^{SUPG}_\tau}$ \\ constant $\tau$\end{tabular} & 
    \begin{tabular}[c]{@{}c@{}}$\mathcal{L^{SUPG}_\tau}$ \\ learnt $\tau$\end{tabular} \\ 
    \hline
    search range & - & $\tau \in [10^{-5}, 10^{-1}]$ & $\tau_g \in [ 5 \cdot10^{-2},10]$ \\ 
    \hline
    optimal parameter & - & $\tau = 10^{-5}$ & $\tau_g = 1$ \\
    \hline
    best $L^2_{\text{err}}$ & $1.693 \cdot 10^{-4}$ & $1.340 \cdot 10^{-4}$ & $1.037 \cdot 10^{-4}$ \\
    \hline
    \end{tabular}
\end{table}
\begin{table}[h]
\centering
\caption{$L^2_{\text{err}}$ for problem $\text{P}_{\text{para}}$ using different loss functionals}
\label{tab:adaptive_tau_accuracy_new_para}
\renewcommand{\arraystretch}{1.2}
    \small
    \begin{tabular}{|l|c|c|c|}
    \hline
    loss & $\mathcal{L}^{\text{hard}}_{\text{var}}$ & 
    \begin{tabular}[c]{@{}c@{}}$\mathcal{L^{SUPG}_\tau}$ \\ constant $\tau$\end{tabular} & 
    \begin{tabular}[c]{@{}c@{}}$\mathcal{L^{SUPG}_\tau}$ \\ learnt $\tau$\end{tabular} \\ 
    \hline
    search range & - & $\tau \in [10^{-4}, 10^{-1}]$ & $\tau_g \in [ 10^{-4},10]$ \\ 
    \hline
    optimal parameter & - & $\tau = 5\cdot 10^{-4}$ & $\tau_g = 1$ \\
    \hline
    best $L^2_{\text{err}}$ & $1.544 \cdot 10^{-4}$ & $1.192 \cdot 10^{-4}$ & $9.043 \cdot 10^{-5}$ \\
    \hline
    \end{tabular}
\end{table}
All errors presented in Tables~\ref{tab:adaptive_tau_accuracy_new_out} and~\ref{tab:adaptive_tau_accuracy_new_para} are considerably smaller than those in Tables~\ref{tab:outflow_layer_params} and~\ref{tab:circular_interior_layer_params}, respectively. For $\text{P}_{\text{out}}$, 
in particular the error in the right lower corner of the domain almost disappeared, 
compare Figures~\ref{fig:Inital_expt} and~\ref{fig:Modified_ansatz_solution_tau}. 
In the case of $\text{P}_{\text{para}}$, the contour lines of the computed solutions are parallel to the $y$-axis now. 
These observations show that a combination of a learnt stabilization parameter and a suitable indicator
function, namely using \eqref{eq: Ansatz_2} and \eqref{eq: Ansatz_3} instead of \eqref{eq: Ansatz_initial}, leads to a tremendous increase of accuracy. 

\subsection{Adaptive Indicator Function}

In our previous set of studies, we used a fixed indicator function to impose hard-constrained boundary conditions 
and we observed in preliminary simulations that the choice of the indicator function might have a 
considerable impact on the accuracy of the computed solutions. To the best of our knowledge, there is no 
proposal available on how to select a priori a good indicator function. 
That's why, we pursue an approach to compute an adaptive indicator function. This function is still of the same principal form as used so far, 
but it is constructed using trainable parameters 
that control the slope of the function near the boundaries. Our approach will enable the neural network to learn not 
only the solution but also the parameters for the indicator functions, resulting in a better fit for the 
specific parameters of the chosen problem.

Concerning problem $\text{P}_{\text{out}}$, we propose an adaptive indicator function of form \eqref{eq: Ansatz_adaptive}, 
more precisely,
\begin{equation}
\begin{split}
     h(x,y) &= \left( 1 - \text{e}^{-\kappa_1 x} \right) \left( 1 - \text{e}^{-\kappa_1 y} \right) \left( 1 - \text{e}^{-\kappa_2 (1 - x)} \right)  \left( 1 - \text{e}^{-\kappa_2 (1 - y)} \right), \\
     \kappa_1 &= 10^{\alpha},\;\kappa_2 = 10^{\beta},
\end{split}
\label{eq: Ansatz_adaptive}
\end{equation}
where $\alpha$ and $\beta$ are learnable parameters controlling the gradient near the boundaries. 
Specifically, $\alpha$ governs the behavior at the inlet boundaries,  while $\beta$ controls the behavior 
at the outflow boundaries. 

Similarly, we introduce an adaptive indicator function for $\text{P}_{\text{para}}$ with three learnable parameters: $\alpha$, $\beta$, and $\gamma$,
namely
\begin{equation}
\begin{split}
     h(x,y) &= \left( 1 - \text{e}^{-\kappa_1 x} \right) \left( 1 - \text{e}^{-\kappa_2 y} \right) \left( 1 - \text{e}^{-\kappa_3 (1 - x)} \right)  \left( 1 - \text{e}^{-\kappa_2 (1 - y)} \right), \\
     \kappa_1 &= 10^{\alpha},\;\kappa_2 = 10^{\beta},\;\kappa_3 = 10^{\gamma}.
\end{split}
\label{eq: Ansatz_adaptive_3}
\end{equation}
In this formulation, $\alpha$ affects the inlet boundary, $\beta$ the characteristic boundaries, and $\gamma$ the outflow boundary. 

The numerical studies with adaptive indicator functions were performed only with the loss functional $\mathcal{L}^{\text{hard}}_{\text{var}}$
from \eqref{eq:vpinn_hard_loss}. Since the steepness of layers depends on the P\'eclet number, we considered different values 
of the diffusion coefficient to investigate how the network's prediction of the indicator 
function depends on the steepness of the layers.

\begin{figure}[t!]
    \centering
\includegraphics[width=1\textwidth]{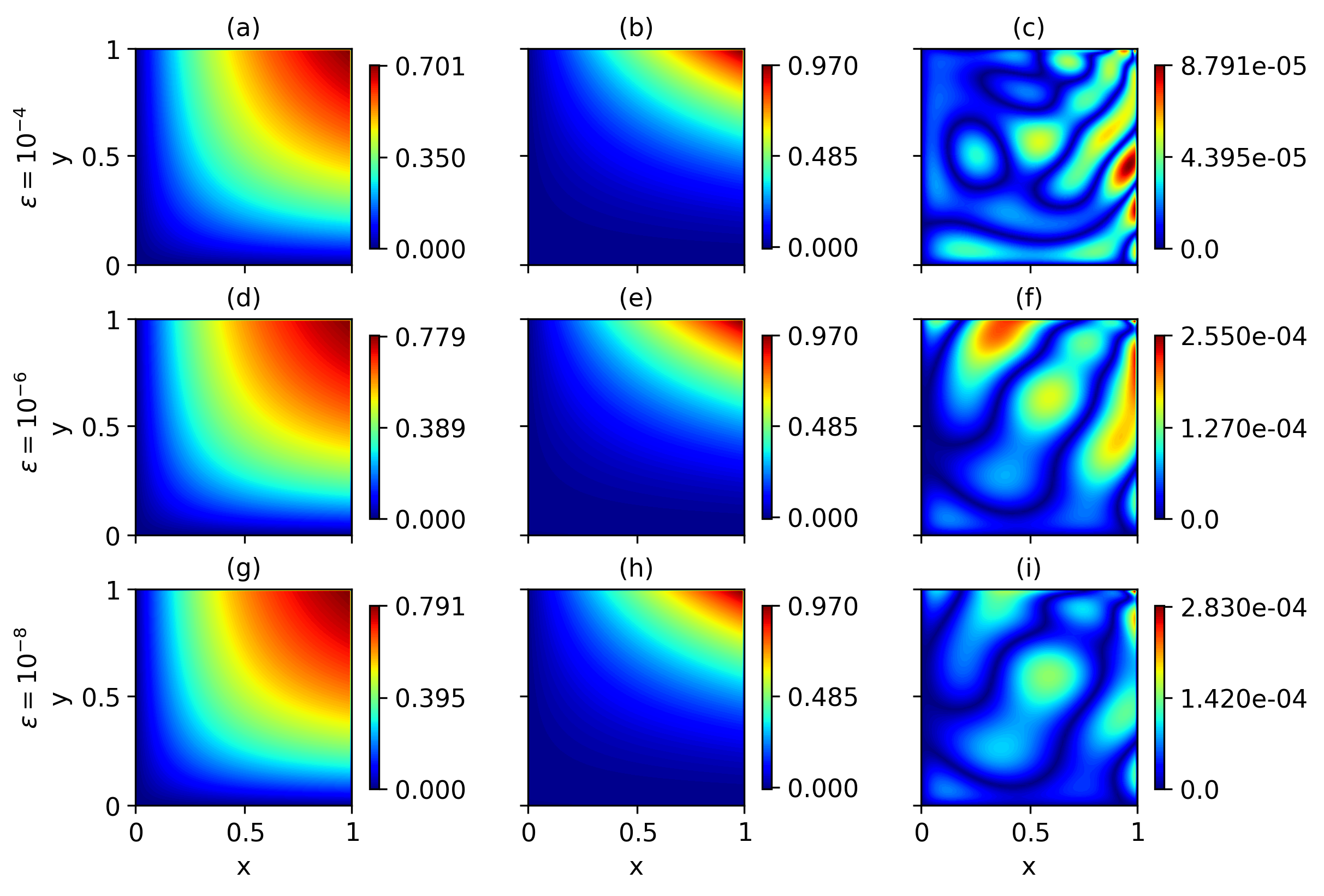}
    \caption{Results for problem $\text{P}_{\text{out}}$ for different $\varepsilon$. The indicator function $h(x,y)$ is shown in (a), (d) and (g), the predicted solution in (b), (e) 
    and (h), and the point-wise error in (c), (f) and (i)} 
    \label{fig:Modified_ansatz_solution_tau_eps}
\end{figure}

\begin{table}[t!]
\caption{Adaptive indicator function parameters and $L^2_{\text{err}}$ for $\text{P}_{\text{out}}$}
\label{tab:adap_indicator_P1}
\centering
\renewcommand{\arraystretch}{1.3}
\small
\begin{tabular}{|c|c|c|c|c|c|c|}
\hline
$\varepsilon$ & \begin{tabular}[c]{@{}c@{}}initial\\ $\alpha$\end{tabular} & \begin{tabular}[c]{@{}c@{}}final\\ $\alpha$\end{tabular} & \begin{tabular}[c]{@{}c@{}}initial\\ $\beta$\end{tabular} & \begin{tabular}[c]{@{}c@{}}final\\ $\beta$\end{tabular}  &  \begin{tabular}[c]{@{}c@{}} initial\\$L^2_{\text{err}}$ \end{tabular}   &  \begin{tabular}[c]{@{}c@{}} final\\$L^2_{\text{err}}$ \end{tabular}\\
\hline
$10^{-4}$ & $1$ & $0.281$ & $2$ & $3.998$ & $2.295\cdot 10^{-1}$  & $4.317\cdot 10^{-5}$\\
\hline
$10^{-6}$ & $1$ & $0.338$ & $3$ & $3.999$  & $1.034\cdot 10^{-1}$ & $5.518\cdot 10^{-5}$\\
\hline
$10^{-8}$ & $1$ & $0.360$ & $4$ & $4.051$  & $9.936\cdot 10^{-5}$ & $5.887\cdot 10^{-5}$\\
\hline
\end{tabular}
\end{table}

Results for $\text{P}_{\text{out}}$ are presented in Figure~\ref{fig:Modified_ansatz_solution_tau_eps} and Table~\ref{tab:adap_indicator_P1}.
The parameters for the indicator function were initialized with $\alpha = 1$ and $\beta = \phi$, where $\phi = -\log_{10}(\varepsilon)/2$. 
For instance, when $\varepsilon = 10^{-8}$, we set $\beta = -\log_{10}(10^{-8})/2 = 4$. It can be seen in Table~\ref{tab:adap_indicator_P1} 
that in fact very different parameters are proposed by  the neural network for the inlet and outflow boundaries, respectively, with the
value at the inlet boundaries being considerably smaller. The parameters do not depend very much on the diffusion coefficient, i.e., on the 
P\'eclet number. The indicator functions have a different shape than the function given in \eqref{eq: Ansatz_2}, which takes values very close to 1
in most parts of $\Omega$. It can be seen that the results obtained with the adapted parameters are (much) more accurate than those with 
the initial parameters. 
Comparing with Table~\ref{tab:adaptive_tau_accuracy_new_out}, one can notice also that the error in the case $\varepsilon=10^{-8}$ is noticeably smaller if the adaptive indicator
function is used. 

\begin{figure}[t!]
    \centering
\includegraphics[width=1\textwidth]{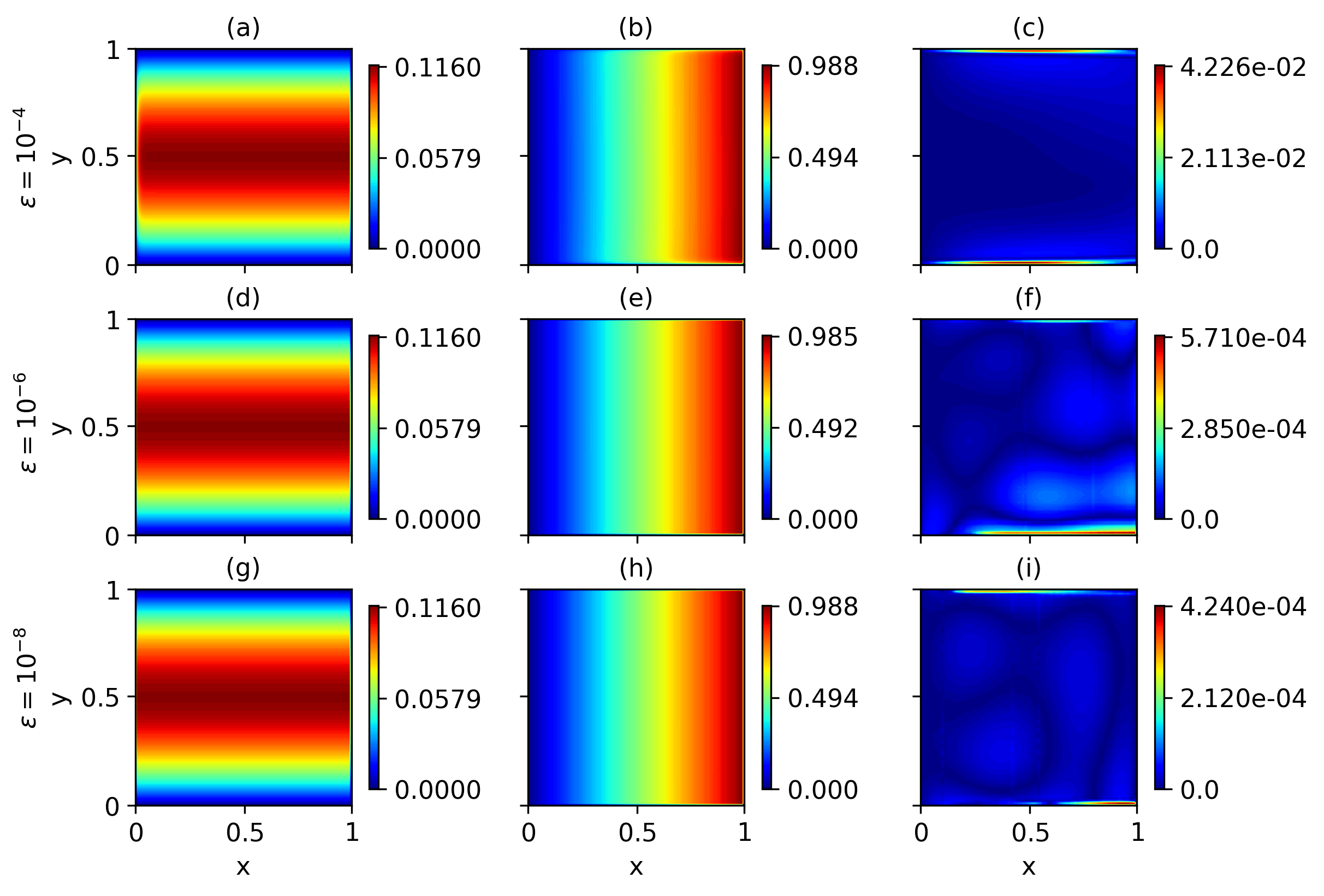}
    \caption{Results for problem $\text{P}_{\text{para}}$ for different $\varepsilon$. The indicator function $h(x,y)$ is shown in (a), (d) and (g), the predicted solution in (b), (e) 
    and (h), and the point-wise error in (c), (f) and (i)} 
    \label{fig:Modified_ansatz_solution_tau_eps_para}
\end{figure}

\begin{table}[t!]
\caption{Adaptive indicator function parameters and $L^2_{\text{err}}$ for $\text{P}_{\text{para}}$}
\label{tab:adap_indicator_P3}
\centering
\renewcommand{\arraystretch}{1.2}
\footnotesize
\begin{tabular}{|c|c|c|c|c|c|c|c|c|c|c|}
\hline
$\varepsilon$ & \begin{tabular}[c]{@{}c@{}}initial\\ $\alpha$\end{tabular} & \begin{tabular}[c]{@{}c@{}}final\\ $\alpha$\end{tabular} & \begin{tabular}[c]{@{}c@{}}initial\\ $\beta$\end{tabular} & \begin{tabular}[c]{@{}c@{}}final\\ $\beta$\end{tabular} & \begin{tabular}[c]{@{}c@{}}initial\\ $\gamma$\end{tabular} & \begin{tabular}[c]{@{}c@{}}final\\ $\gamma$\end{tabular} &  \begin{tabular}[c]{@{}c@{}}initial\\$L^2_{\text{err}}$\end{tabular} & \begin{tabular}[c]{@{}c@{}}final \\$L^2_{\text{err}}$\end{tabular} \\
\hline
$10^{-4}$ & 1 & 2.129 & 0 & -0.086 & 2 & 3.927  & $5.124\cdot10^{-1}$& $1.075\cdot10^{-2}$\\
\hline
$10^{-6}$ & 1 & 3.295 & 0 & -0.088 & 3 & 3.921  & $2.706\cdot10^{-1}$& $1.066\cdot10^{-4}$\\
\hline
$10^{-8}$ & 2 & 3.380 & 0 & -0.079 & 4 & 4.004  & $1.534\cdot10^{-2}$& $6.007\cdot10^{-5}$\\
\hline
\end{tabular}
\end{table}

For problem $\text{P}_{\text{para}}$, whose solution possesses different types of layers and whose ansatz for the 
indicator function is given in \eqref{eq: Ansatz_adaptive_3}, the obtained results are
given in Figure~\ref{fig:Modified_ansatz_solution_tau_eps_para} and Table~\ref{tab:adap_indicator_P3}. The evaluation 
comes to similar conclusions as for the other example. Again, the learnt parameters for the indicator function 
depend only slightly on the P\'eclet number. The visualizations of the indicator functions look very much alike. 
Comparing with Table~\ref{tab:adaptive_tau_accuracy_new_para}, one can observe that the results for $\varepsilon=10^{-8}$ is considerably more accurate 
for the approach with adaptive indicator function.

\section{Conclusion}\label{sec:conclusion}

FastVPINNs are a very efficient approach for computing approximations of
solutions of boundary value problems. In this paper, this approach was 
combined with two methods for improving the accuracy of computed solutions 
of convection-diffusion-reaction problems in the convection-dominated regime. 
First, the residual loss was augmented with a SUPG stabilization term and an 
architecture was proposed that computes a spatially varying stabilization 
parameter. And second, an architecture was developed that learns appropriate 
parameters for defining a good indicator function for imposing hard-constrained
Dirichlet boundary conditions. Both extensions led to a considerable increase 
of the accuracy of the computed solutions. In contrast, the inclusion of a 
regularization term with the weights of the network was not helpful for obtaining 
solutions with higher accuracy.


\newpage



\bibliography{ref}
\end{document}